\documentclass[11pt]{amsart}
\usepackage{amsmath,amssymb,amsthm,mathtools}
\usepackage[colorlinks=true,linkcolor=black,citecolor=black,urlcolor=blue]{hyperref}
\usepackage[margin=1.07in]{geometry}
\newtheorem{theorem}{Theorem}[section]
\newtheorem{lemma}[theorem]{Lemma}
\newtheorem{corollary}[theorem]{Corollary}

\newcommand{\supp}{\operatorname{supp}}
\newcommand{\one}{\mathbf 1}

\title[The ball fixed point property for continuous functions on locally compact spaces]{The ball fixed point property for continuous functions on locally compact spaces}
\author{Lily DiPaulo}
\date{October 30, 2026}
\begin{document}
\begin{abstract}
Let $K$ be a compact Hausdorff space and $p\in K$.
We prove that the Banach space
$C(K,p)=\{f\in C(K,\mathbb R):f(p)=0\}$
fails the ball fixed point property (BFPP) whenever
$p$ is non-isolated, answering a question of
Avil\'es, Jap\'on, Lennard, Mart\'inez-Cervantes
and Stawski. Consequently, $C_0(L)$ fails the
BFPP for every noncompact locally compact Hausdorff
space $L$. Together with Barroso's characterization
for $C(K)$, this yields complete topological
classifications of the BFPP for $C_0(L)$ and $C(K,p)$.
We also resolve a second question of Avil\'es
et al., showing that the BFPP of $C(K)$ implies
that of $C(K,p)$ for every compact Hausdorff
space $K$ and every $P$-point $p\in K$ if
and only if there are no measurable cardinals.
\end{abstract}
\maketitle

\section{Introduction}
A Banach space $X$ has the \emph{ball fixed point property} (BFPP) if every
nonexpansive mapping from the closed unit ball of $X$ into itself has a fixed point. Unless otherwise specified, $C(K)$ and $C_0(L)$ denote spaces of real-valued continuous functions equipped with the supremum norm.

For a compact Hausdorff space $K$ and $p\in K$, define
\[
C(K,p)=\{f\in C(K):f(p)=0\}.
\]
Recall that $p\in K$ is a $P$-point if every
$f\in C(K,p)$ vanishes on a neighborhood of $p$. Avil\'es, Jap\'on, Lennard, Mart\'inez-Cervantes and Stawski
\cite[Section~6]{avil} showed that $C(K,p)$ fails BFPP if $p$ is not a
$P$-point, and asked whether the same conclusion holds for every
non-isolated $p$. Indeed, we show that $C(K,p)$ fails the BFPP for every non-isolated point $p$, regardless of whether $p$ is a $P$-point. We also show that the converse of
\cite[Theorem~6.1]{avil}, restricted to spaces
$C(K)$ with the BFPP, holds universally if
and only if there are no measurable cardinals.

Our construction adapts the transfinite propagation technique introduced by Barroso \cite{barroso}. The condition $f(p)=0$ allows us to replace the two profiles used in that work with a single decreasing profile. We also use the positive linear contraction of \cite[Lemma~3.2]{barroso}.

\section{The ball fixed point property in $C(K,p)$}
\begin{theorem}\label{thm:main}
Let $K$ be a compact Hausdorff space, and suppose that $p\in K$ is non-isolated.
Then $C(K,p)$ fails the BFPP.
\end{theorem}

\begin{lemma}\label{lem:family}
Let $K$ be a compact Hausdorff space and
$p\in K$ non-isolated. Then there are a limit
ordinal $\tau$ and an increasing family
$(u_\alpha)_{\alpha<\tau}$ in $C(K,[0,1])$ such that
$u_0=0$, each $u_\alpha$ vanishes near $p$, and, with
$F_\alpha=\{u_\alpha=1\}$,
\begin{equation}\label{eq:closure}
p\in\overline{\bigcup_{\alpha<\tau}F_\alpha}.
\end{equation}
\end{lemma}
\begin{proof}
Enumerate $K\setminus\{p\}=\{x_\alpha:\alpha<\kappa\}$,
where $\kappa$ is an infinite cardinal, and set $u_0=0$.
We proceed by transfinite recursion. Having constructed
$(u_\beta)_{\beta<\alpha}$, define
\[
H_\alpha=\overline{\bigcup_{\beta<\alpha}\supp u_\beta}.
\]
Stop if $p\in H_\alpha$. Otherwise, for $0<\alpha<\kappa$,
normality of $K$ yields $u_\alpha\in C(K,[0,1])$
such that $u_\alpha=1$ on $H_\alpha$, $u_\alpha$
vanishes on a neighborhood of $p$, and
$u_\alpha(x_\beta)=1$ whenever $\alpha=\beta+1$.
By construction, the functions are increasing and
$H_{\alpha+1}=\supp u_\alpha$.

The recursion must stop by stage $\kappa$. Otherwise, we have that
\[
x_\beta\in F_{\beta+1}\subseteq H_\kappa
\qquad(\beta<\kappa).
\]
Since $H_\kappa$ is closed and $p$ is non-isolated,
this implies $p\in H_\kappa$, a contradiction.

Let $\tau$ be the first stopping ordinal. Since
$H_{\alpha+1}=\supp u_\alpha$ does not contain $p$,
the ordinal $\tau$ is a limit. Moreover,
\[
\supp u_\alpha\subseteq F_{\alpha+1}
\qquad(\alpha<\tau).
\]
Consequently,
\[
p\in H_\tau
=\overline{\bigcup_{\alpha<\tau}\supp u_\alpha}
\subseteq
\overline{\bigcup_{\alpha<\tau}F_\alpha}.
\]
\end{proof}

\begin{lemma}\label{lem:synthesis}
Let $K$ be a compact Hausdorff space, let $\tau$
be a limit ordinal, and let
$(u_\alpha)_{\alpha<\tau}$ be an increasing
family in $C(K,[0,1])$. There is a positive linear contraction
\[
S:C_0([0,\tau))\longrightarrow C(K),\qquad
S\one_{[0,\alpha]}=u_\alpha \quad (\alpha<\tau).
\]
If $p\in K$ and $u_\alpha(p)=0$ for every
$\alpha<\tau$, then the range of $S$
is contained in $C(K,p)$.
\end{lemma}
\begin{proof}
Let $e_\alpha=\one_{[0,\alpha]}$. These functions
are linearly independent, and their span
$D=\operatorname{span}\{e_\alpha:\alpha<\tau\}$
is dense in $C_0([0,\tau))$ by the
Stone--Weierstrass theorem. Define
$S:D\to C(K)$ by $Se_\alpha=u_\alpha$.

For $\alpha_1<\cdots<\alpha_n$, the functions
\[
e_{\alpha_1},\quad
e_{\alpha_2}-e_{\alpha_1},\quad\ldots,\quad
e_{\alpha_n}-e_{\alpha_{n-1}}
\]
are nonnegative and have pairwise disjoint supports.
Their images under $S$ are nonnegative, so $S$ is
positive. Moreover, every
$g\in\operatorname{span}\{e_{\alpha_1},\ldots,e_{\alpha_n}\}$
satisfies
\[
-\|g\|_\infty e_{\alpha_n}
\le g\le
\|g\|_\infty e_{\alpha_n}.
\]
By positivity,
\[
|Sg|\le\|g\|_\infty u_{\alpha_n}
\le\|g\|_\infty.
\]
Therefore $S$ extends uniquely to a positive linear
contraction on $C_0([0,\tau))$.

If $u_\alpha(p)=0$ for every $\alpha<\tau$,
then $Sg(p)=0$ for every $g\in D$, and hence
for every $g\in C_0([0,\tau))$ by continuity.
\end{proof}

\begin{proof}[Proof of Theorem~\ref{thm:main}]
Let $(u_\alpha)_{\alpha<\tau}$ and $S$ be as in
Lemmas~\ref{lem:family} and \ref{lem:synthesis}, and define
\[
\mathcal C=\{f\in C(K,p):0\le f\le1\}.
\]
For $f\in\mathcal C$ and $\alpha<\tau$, let
\[
m_\alpha(f)=
\begin{cases}
\min_{F_\alpha}f,&F_\alpha\ne\emptyset,\\
1,&F_\alpha=\emptyset.
\end{cases}
\]
Since the sets $F_\alpha$ are increasing and
$p\in\overline{\bigcup_{\alpha<\tau}F_\alpha}$,
continuity of $f$ at $p$ gives
\begin{equation}\label{eq:profile}
m_0(f)=1,\qquad
m_\alpha(f)\ge m_\beta(f)\quad(\alpha\le\beta),
\qquad
\inf_{\alpha<\tau}m_\alpha(f)=0.
\end{equation}

Define $h_f:[0,\tau)\to[0,1]$ by
\begin{equation}\label{eq:shift}
\begin{aligned}
h_f(0)&=1,\\
h_f(\alpha+1)&=m_\alpha(f),\\
h_f(\lambda)&=\inf_{\alpha<\lambda}m_\alpha(f)
\end{aligned}
\end{equation}
for nonzero limit ordinals $\lambda<\tau$.
By \eqref{eq:profile}, $h_f$ is decreasing,
continuous in the order topology, and vanishes at
infinity. Moreover, since minima over fixed compact
sets and infima of bounded families are nonexpansive, we have that
\begin{equation}\label{eq:lip}
\|h_f-h_g\|_\infty\le\|f-g\|_\infty
\qquad(f,g\in\mathcal C).
\end{equation}

Define $T:\mathcal C\to\mathcal C$ by $Tf=Sh_f$.
Positivity and contractivity of $S$, together with
\eqref{eq:lip}, imply that $T$ is nonexpansive.
Writing $e_\alpha=\one_{[0,\alpha]}$, we have
\[
h_f\ge h_f(\alpha)e_\alpha,
\qquad
Tf\ge h_f(\alpha)u_\alpha,
\]
and consequently
\begin{equation}\label{eq:ineq}
m_\alpha(Tf)\ge h_f(\alpha)
\qquad(\alpha<\tau).
\end{equation}

Suppose $Tf=f$. By \eqref{eq:shift} and
\eqref{eq:ineq},
\[
m_{\alpha+1}(f)\ge m_\alpha(f)
\qquad(\alpha+1<\tau),
\]
and
\[
m_\lambda(f)\ge
\inf_{\alpha<\lambda}m_\alpha(f)
\]
for every nonzero limit ordinal $\lambda<\tau$.
Since $m_0(f)=1$, transfinite induction gives
\[
m_\alpha(f)=1\qquad(\alpha<\tau),
\]
contradicting \eqref{eq:profile}. Hence $T$ has no
fixed point.

Finally, define $R:B_{C(K,p)}\to\mathcal C$ such that
\[
Rf=f^+=\max\{f,0\}.
\]
Then $R$ is a nonexpansive retraction. So $T\circ R$
is a fixed-point-free nonexpansive mapping on
$B_{C(K,p)}$ and consequently $C(K,p)$ fails the BFPP.
\end{proof}

\section{Corollaries}
\begin{corollary}\label{cor:c0}
If $L$ is a noncompact locally compact Hausdorff space, then $C_0(L)$
fails the BFPP.
\end{corollary}
\begin{proof}
Let $K=L\cup\{\infty\}$ be the one-point compactification of $L$.
Since $\infty$ is
non-isolated and $C_0(L)$ is isometric to
$C(K,\infty)$, Theorem~\ref{thm:main} implies
that $C_0(L)$ fails the BFPP.
\end{proof}

A compact space is called \emph{extremally disconnected} if the closure
of every open subset is open. Barroso \cite[Corollary~5.4]{barroso}
proved that $C(K)$ has BFPP precisely for extremally disconnected
compact spaces $K$. Together with Corollary~\ref{cor:c0}, this yields
the following classification.

\begin{corollary}\label{cor:classification}
For a locally compact Hausdorff space $L$,
\[
C_0(L)\text{ has BFPP}
\quad\Longleftrightarrow\quad
L\text{ is compact and extremally disconnected}.
\]
For compact Hausdorff $K$ and $p\in K$, equivalently,
\[
C(K,p)\text{ has BFPP}
\quad\Longleftrightarrow\quad
p\text{ is isolated and }K\text{ is extremally disconnected}.
\]
\end{corollary}
\begin{proof}
The characterization of $C_0(L)$ follows from
Corollary~\ref{cor:c0} and \cite{barroso}.
For $C(K,p)$, observe that
$C(K,p)$ is isometric to $C_0(K\setminus\{p\})$.
The space $K\setminus\{p\}$ is compact precisely
when $p$ is isolated, in which case it is
extremally disconnected precisely when $K$ is extremally disconnected.
\end{proof}

\begin{corollary}\label{cor:ppoints}
The following are equivalent:
\begin{enumerate}
    \item There are no measurable cardinals.
    \item For every compact Hausdorff space $K$ and every
    $P$-point $p\in K$, if $C(K)$ has the BFPP, then
    $C(K,p)$ has the BFPP.
\end{enumerate}
\end{corollary}

\begin{proof}
Suppose there are no measurable cardinals. If $C(K)$
has the BFPP, then $K$ is extremally disconnected
by \cite{barroso}, and every $P$-point of $K$ is
isolated by \cite[Section~4.4]{terasawa}.
The conclusion follows from
Corollary~\ref{cor:classification}.

Conversely, suppose that a measurable cardinal
$\kappa$ exists. Let $p\in\beta\kappa$ correspond
to a nonprincipal countably complete ultrafilter
on $\kappa$. Then $p$ is a non-isolated $P$-point
of $\beta\kappa$. Since $\beta\kappa$ is extremally
disconnected, $C(\beta\kappa)$ has the BFPP,
whereas $C(\beta\kappa,p)$ fails the BFPP by
Theorem~\ref{thm:main}.
\end{proof}

\begin{corollary}\label{cor:ideals}
Let $I$ be a closed ideal of $C(K)$. Then $I$
has the BFPP if and only if
\[
I=\one_U C(K)
\]
for some clopen, extremally disconnected
subset $U$ of $K$.
\end{corollary}
\begin{proof}
Every closed ideal $I$ of $C(K)$ has the form
\[
I=\{f\in C(K):f|_F=0\}\cong C_0(K\setminus F)
\]
for some closed subset $F\subseteq K$. Let $U=K\setminus F$.
By Corollary~\ref{cor:classification}, $I$ has the BFPP
if and only if $U$ is compact and extremally disconnected.
Since $U$ is open in $K$, it is compact if and only if
it is clopen, in which case $I=\one_U C(K)$.
\end{proof}

\begin{corollary}\label{cor:vector}
Let $L$ be a noncompact locally compact Hausdorff space
and $E$ a nonzero real or complex Banach space. Then
the underlying real Banach space $C_0(L,E)$ fails
the BFPP.
\end{corollary}

\begin{proof}
Choose $e\in E$ with $\|e\|=1$ and a norm-one
real-linear functional $\varphi:E\to\mathbb R$
satisfying $\varphi(e)=1$. Define
\[
\begin{aligned}
P:C_0(L,E)&\to C_0(L),
& (Pf)(t)&=\varphi(f(t)),\\
J:C_0(L)&\to C_0(L,E),
& (Jg)(t)&=g(t)e.
\end{aligned}
\]
Both operators are real-linear contractions, and
$PJ=\operatorname{id}_{C_0(L)}$.

By Corollary~\ref{cor:c0}, there exists a
fixed-point-free nonexpansive mapping
$A:B_{C_0(L)}\to B_{C_0(L)}$.
Consequently,
\[
JAP:B_{C_0(L,E)}\to B_{C_0(L,E)}
\]
is nonexpansive. If $f=JAPf$, then
\[
Pf=PJAPf=APf,
\]
contradicting the absence of fixed points of $A$.
\end{proof}

\begin{corollary}\label{cor:cstar}
Every nonzero nonunital commutative complex $C^*$-algebra
fails the BFPP when regarded as a real Banach space.
\end{corollary}

\begin{proof}
By the commutative Gelfand--Naimark theorem, such an
algebra is isometrically isomorphic to $C_0(L,\mathbb C)$
for some locally compact Hausdorff space $L$.
Since the algebra is nonunital, $L$ is noncompact,
and the conclusion follows from
Corollary~\ref{cor:vector}.
\end{proof}

\section*{Acknowledgements}

The author used OpenAI's ChatGPT 6.1 Sol
to assist with the development and review of
proofs and the preparation of this manuscript.
The author assumes full responsibility for
the results and final content.

\end{document}